\documentclass[12pt]{article}
\usepackage{amssymb,amsmath,amsfonts,amsthm,amsxtra}
\usepackage{bm}
\usepackage{cite}
\usepackage{array}
\usepackage{graphicx}

\def\e{\mathrm e}                                               %+
\def\td{\tilde d}                                                %+

\def\ze{\zeta}                                                  %+
\def\RR{\mathop{\cal R}\nolimits}                               %+
\def\CC{\mathop{\cal C}\nolimits}                               %+
\def\cdc{,\ldots,}                                              %+
\def\1n{1\cdc n}                                                %+
\def\eq#1{\begin{equation}#1\end{equation}}                     %+
\def\eqs*#1{\begin{eqnarray*}#1\end{eqnarray*}}                 %+
\def\eqss#1{\begin{eqnarray}#1\end{eqnarray}}                   %+
\newtheorem{thm}{Theorem}{\bfseries}{\itshape}                  %+
\newtheorem{corol}{Corollary}{\bfseries}{\itshape}              %+
\newtheorem{lemma}{Lemma}{\bfseries}{\itshape}                  %+
{\bfseries}{\itshape}              %+
\newtheorem{defin}{Definition}{\bfseries}{\upshape}             %+
\newtheorem{remark}{Remark}{\bfseries}{\upshape}                %+
\newtheorem{example}{Example}{\bfseries}{\upshape}              %+
\def\proof{{\noindent\bf Proof. }}                              %+
\def\R{{\mathbb R}}                                             %+
\def\ovec{\overrightarrow}                                      %+
\def\tot{\rightleftarrows}                                      %+
\def\T{{\xz\rm\scriptscriptstyle T}\xz}                         %+
\def\W{{\hspace{.05em}\rm\scriptscriptstyle W}}                 %+

\def\SS{{\cal S}}                                               %+
\def\xy{\hspace{.07em}}                                         %+
\def\xz{\hspace{-.07em}}                                        %+
\def\ms{{\mathstrut}}                                           %+
\def\tr{\operatorname{tr}}                                  %+
\def\_#1{{^{}_{#1}}}                     %+
\def\j{{\bar\jmath}}                                            %+
\def\i{{\bar\imath}}                                            %+
\def\jj{{\bar\jmath\bar\jmath}}                                 %+
\def\ii{{\bar\imath\hspace*{ .075ex}\bar\imath}}                 %+
\def\ij{{\bar\imath\hspace*{-.065ex}\bar\jmath}}                 %+
\def\ji{{\bar\jmath\hspace*{ .130ex}\bar\imath}}                 %+
\def\Up#1{\vspace*{-#1em}}                                       %+
\def\AA{{\cal A}}                                               %+
\def\G{\Gamma}                                                  %+
\def\sm{\xz\smallsetminus\xz}                                   %+
\def\mi{\xz-\xz}                                                %+
\def\w{{\mathrm w}}                                             %+
\def\too{\!\to\!}                                               %+
\def\rr{{\mathrm r}}                                            %+
\def\ai{_\ms^\ms}                                               %+
\def\ak{_\ms}                                                   %+
\def\aj{^\ms}                                                   %+

\author{Pavel Chebotarev\footnotemark[1]
        \and Michel Deza\footnotemark[2]
}

\title{\vspace{-.8em}
A topological interpretation of the walk distances}
\date{}

\def\thefootnote{\fnsymbol{footnote}}
\begin{document}
\footnotetext[1]{{Institute of Control Sciences of the Russian Academy of Sciences,}
                               {65 Profsoyuznaya Street, Moscow 117997, Russia,}
                               E-mail:\:{\footnotesize\tt chv@member.ams.org.}}
\footnotetext[2]{{Laboratoire de Geometrie Appliquee, LIGA, Ecole Normale Superieure,}
                               {45, rue d'Ulm, F-75230, Paris, Cedex 05, France,}
                               E-mail:\;{\footnotesize\tt Michel.Deza@ens.fr.}}
\maketitle
\def\thefootnote{\arabic{footnote}}
%\vspace{-1em}
\begin{abstract}
The walk distances in graphs have no direct interpretation in terms of walk weights, since they are introduced via the \emph{logarithms\/} of walk weights. Only in the limiting cases where the logarithms vanish such representations follow straightforwardly. The interpretation proposed in this paper rests on the identity ${\ln\det B=\tr\ln B}$ applied to the cofactors of the matrix $I-tA,$ where $A$ is the weighted adjacency matrix of a weighted multigraph and $t$ is a sufficiently small positive parameter. In addition, this interpretation is based on the power series expansion of the logarithm of a matrix. Kasteleyn~\cite{Kasteleyn67} was probably the first to apply the foregoing approach %aforementioned identities
to expanding the determinant of $I-A.$ We show that using a certain linear transformation the same approach can be extended to the cofactors of $I-tA,$ which provides a topological interpretation of the walk distances.


\medskip
%\bigskip
%15A24, %Matrix equations and identities
%15B48  %Positive matrices and their generalizations; cones of matrices (before: 15A48)
%15B51%, %Stochastic matrices (before: 15A51)
\end{abstract}

\section{Introduction}

The walk distances for graph vertices were proposed in~\cite{Che11AAM} and studied in~\cite{Che12DAM}.
%They form a class of distances which
Along with their modifications they generalize~\cite{Che12DAM} the logarithmic forest distances \cite{Che11DAM}, resistance distance, shortest path distance, and the weighted shortest path distance.
The walk distances are graph-geodetic: for a distance\footnote{In this paper, a \emph{distance\/} is assumed to satisfy the axioms of metric.} $d(i,j)$ in a graph $G$\/ this means that $d(i,j)+d(j,k)=d(i,k)$ if and only if every path in $G$ connecting $i$ and $k$ visits~$j.$

%Two classical distances for graph vertices are the shortest path distance~\cite{BuckleyHarary90} and the resistance distance~\cite{GvishianiGurvich87En,KleinRandic93}, which is proportional to the commute time distance~\cite{GobelJagers74RandomW}. Recently, a need for a wider variety of graph distances has been strongly felt (see, e.g.,~\cite{YenSaerensShimbo08,DezaDeza09,LuxburgRadlHein09,Tang10PhD} among many others).%\nocite{Higham08,Che12DAM,DezaDeza09,Che11DAM}

%In \cite{CheSha98} we applied the inverse covariance mapping to the matrices of the walk weights $\sum_{k=0}^\infty (tA)^k,$ where $A$ is the adjacency matrix of a graph, and showed that this leads to distances whenever the positive parameter $t$ is sufficiently small. However, these distances are not graph-geodetic.

%It is well known that the resistance distance coincides with the shortest path distance on any tree. This implies that the 
It is well known that the resistance distance between two adjacent vertices in a tree is equal to~$1.$ In contrast to this, the walk distances take into account the centrality of vertices. For example, any walk distance between two central adjacent vertices in a path turns out \cite{Che12DAM} to be less than that between two peripheral adjacent vertices. This property may be desirable in some applications including machine learning, mathematical chemistry, the analysis of social and biological networks, etc.

In the present paper, we obtain a topological interpretation of the simplest walk distances.
Such an interpretation is not immediate from the definition, since the walk distances are introduced via the \emph{logarithms\/} of walk weights. Only in the limiting cases where the logarithms vanish such representations follow straightforwardly~\cite{Che12DAM}. The interpretation we propose rests on the identity ${\ln\det B=\tr\ln B}$ applied to the cofactors of the matrix $I-tA,$ where $A$ is the weighted adjacency matrix of a weighted multigraph and $t$ is a sufficiently small positive parameter. In addition, it is based on the power series expansion of the logarithm of a matrix. We do not employ these identities explicitly; instead, we make use of a remarkable result by Kasteleyn~\cite{Kasteleyn67} based on them. More specifically, Kasteleyn obtained an expansion of the determinant of $I-A$ and the logarithm of this determinant. We show that using a certain linear transformation the same approach can be extended to the cofactors of $I-tA,$ which provides a topological interpretation of the walk distances.

\section{Notation}
\label{s_notat}

In the graph definitions we mainly follow~\cite{Harary69}.
Let $G$ be a weighted multigraph (a weighted graph where multiple edges are allowed) with vertex set $V(G)=V,$ $|V|=n>2,$ %$V(G)=\{\1n\},$  $n>1$
and edge set~$E(G)$. Loops are allowed; we assume that $G$ is connected. For brevity, we will call $G$ a \emph{graph}.
For ${i,j\in V(G),}$ let $n_{ij}\in\{0,1,\ldots\}$ be the number of edges incident to both $i$ and $j$ in~$G$; for every ${q\in\{\1n_{ij}\}}$, $w_{ij}^q>0$ is the weight of the $q\/$th edge of this type. Let
\eq{
\label{e_aij}%+
a_{ij}=\sum_{q=1}^{n_{ij}}w_{ij}^q
}
(if $n_{ij}=0,$ we set $a_{ij}=0$) and $A=(a_{ij})_{n\times n};$ $A$~is the symmetric \emph{weighted adjacency matrix\/} of~$G$. %\emph{matrix of total edge weights} %for all pairs of vertices.
In what follows, all matrix entries are indexed by the vertices of~$G.$ This remark is essential when submatrices are considered: say, ``the $i$th column'' of a submatrix of~$A$ means ``the column corresponding to the vertex $i$ of~$G$'' rather than just the ``column number~$i$.''

%\smallskip
By the \emph{weight\/} of a graph $G$, $w(G)$, we mean the product of the weights of all its edges. If $G$ has no edges, then $w(G)=1$.
The weight of a set $\SS$ of graphs, $w(\SS)$, is the total weight (the sum of the weights) of its elements; %the graphs belonging to~$\SS$; the weight of the empty set is zero.
$w(\varnothing)=0$. %If the weights of all edges are unity, i.e.\ the graphs in $\SS$ are actually unweighted, then $w(\SS)$ reduces to the cardinality of~$\SS$.
%-%The weights of sequences of vertices and edges and of the sets of sequences are defined similarly.

%\smallskip
%For $v_0,v_m\in V(G),$ a $v_0\to v_m$ \emph{path\/} %(\emph{simple path})
%in $G$ is an alternating sequence of vertices and edges $v_0,\e_1,v_1\cdc\e_m,v_m$ where all vertices are distinct and each $\e_i$ is a $(v_{i-1},v_i)$ edge. The unique $v_0\to v_0$ path is the ``sequence''\,$v_0$ having no edges.
%The \emph{length\/} of a path is the number $m$ of its edges. The \emph{weight\/} of a path is the product of the weights of its edges.
%The weight of a $v_0\to v_0$ path is~1.
%
%Similarly,
For $v_0,v_m\in V(G),$ a $v_0\to v_m$ \emph{walk\/} in $G$ % (also called a \emph{route})
is an arbitrary alternating sequence of vertices and edges $v_0,\e_1,v_1\cdc\e_m,v_m$ where each $\e_i$ is a $(v_{i-1},v_i)$ edge. %Some (or all) edges in a walk can be loops.
The \emph{length\/} of a walk is the number $m$ of its edges (including loops and repeated edges). The \emph{weight\/} of a walk is the product of the $m$ weights of its edges. The weight of a set of walks is the total weight of its elements. By definition, for any vertex $v_0$, there is one $v_0\to v_0$ walk $v_0$ with length $0$ and weight~1.
%\smallskip

We will need some special types of walks. A~\emph{hitting $v_0\to v_m$ walk\/} is a $v_0\to v_m$ walk containing only one occurrence of~$v_m.$
A~$v_0\to v_m$ walk is called \emph{closed\/} if $v_m=v_0$ and \emph{open\/} otherwise. The \emph{multiplicity\/} of a closed walk is the maximum $\mu$ such that the walk is a $\mu$-fold repetition of some walk.

%Two closed walks are \emph{similar\/} if one of them can be obtained from the other by a cyclic shift.
We say that two closed walks of non-zero length are \emph{phase twins\/} if the edge sequence $\e_1,\e_2,\cdc\e_m$ of the first walk can be obtained from the edge sequence $\e'_1,\e'_2,\cdc\e'_m$ of the second one by a cyclic shift.
For example, the walks $v_0,\e_1,v_1,\e_2,v_2,\e_3,v_0$ and $v_2,\e_3,v_0,\e_1,v_1,$ $\e_2,v_2$ are phase twins. A~\emph{circuit\/} \cite{Kasteleyn67,HararySchwenk79}
in $G$ is any equivalence class of phase twins. The \emph{mul\-tiplicity\/} of a circuit is the multiplicity of any closed walk it contains (all such walks obviously have the same multiplicity).
A~walk (circuit) whose multiplicity exceeds $1$ is \emph{periodic}.

%A~$v_0\to v_0$ cycle is called a \emph{$v_0\tot v_m$ commute cycle\/} if it contains $v_m$ and has no occurrences of $v_0$ strictly between the first appearance of $v_m$ and the final appearance of~$v_0.$

Let $r_{ij}$ be the weight of the set $\RR^{ij}$ of all $i\too j$ walks in $G$ provided that this weight is finite. $R=R(G)=(r_{ij})_{n\times n}\!\in\xz\R^{n\times n}$ will be referred to as the \emph{matrix of the walk weights\/} of~$G$.

It was shown in \cite{Che11AAM} that if $R$ exists then it \emph{determines a transitional measure in $G$}, that is, (i) it satisfies the transition inequality
\eq{
\label{e_te}
r_{ij}\,r_{\!jk}\le r_{ik}\,r_{\!jj},\quad i,j,k=1\cdc n
}
and (ii) $r_{ij}\,r_{\!jk}=r_{ik}\,r_{\!jj}\xy$ if and only if every path from $i$ to $k$ visits~$j.$

%%\smallskip
%By $d^\s(i,j)$ we denote the \emph{shortest path distance}, i.e., the length of a shortest path between $i$ and $j$ in~$G.$
%The \emph{weighted shortest path distance\/} $d^\ws(i,j)$ is defined as follows:\footnote{This formula corrects Eq.~(6.2) in~\cite{KleinRandic93}; cf.\ \cite[Section\:4]{Klein10JMC}.} %; cf.\ the first inequality in \cite[p.~261]{DezaDeza09}.}
%\eqss{
%\label{e_dws}%+
%d^\ws(i,j)=\min_\pi\sum_{\e\xy\in\xy E(\pi)}l_\e,
%}
%where the minimum is taken over all paths $\pi$ from $i$ to $j$ and the sum is over all edges $\e$ in~$\pi$; $l_\e=1/w_\e$ is sometimes called the \emph{weighted length\/} of the edge $\e$, where $w_\e$ is the weight of this edge (see, e.g.,~\cite{Cinkir11EJC}). In the theory of electrical networks, $l_\e$ is identified with the \emph{resistance\/} of the edge~$\e.$

\def\baselinestretch{1.0}
\section{The walk distances}
\label{s_prel}

For any $t>0,$ consider the graph $tG$ obtained from $G$ by multiplying all edge weights by~$t.$
%By $tG$ we denote the graph obtained from $G$ by multiplying all edge weights by $t$ $(t>0).$
If the matrix of the walk weights of $tG,\,$ $R_t=R(tG)=(r_{ij}(t))_{n\times n},$ exists, then\footnote{In the more general case of weighted \emph{digraphs}, the $ij$-entry of the matrix $R_t-I$ is called the \emph{Katz similarity\/} between vertices $i$ and~$j$. Katz~\cite{Katz53} proposed it to evaluate the social status taking into account all $i\too j$ paths. Among many other papers, this index was studied in \cite{Thompson58,Kasteleyn67}.}
%\footnote{For an early study of the graph proximity measure $\sum_{k=0}^\infty(tA)^k,$ we refer the reader to \cite{Katz53,Thompson58,Ponstein66,Kasteleyn67,Taylor}. More recently, it was explored in \cite{CheSha98,YenSaerensShimbo08,ChelnokovZefirova09,EstradaHigham10}.}
\eq{
\label{e_Rt}%+
R_t=\sum_{k=0}^\infty(tA)^k=(I-tA)^{-1},
}
where $I$ denotes the identity matrix of appropriate dimension.

By assumption, $G$ is connected, while its edge weights are positive, so $R_t$ is also positive. Apply the logarithmic transformation to the entries of~$R_t,$ namely, consider the matrix
\eq{
\label{e_Ha}%+
H_t=\ovec{\ln R_t},
}
where
$\ovec{\varphi(S)}$ stands for elementwise operations, i.e., operations applied to each entry of a matrix $S$ separately. Finally, consider the matrix
\eq{
\label{e_Da}%+
D_t=\frac{1}{2}\xy(h_t{\bm1}^\T+\bm1 h_t^\T-H_t-H_t^\T), %)-H_t,
}
where $h_t$ is the column vector containing the diagonal entries of $H_t$, ${\bm1}$ is the vector of ones of appropriate dimension, and %$H_t^\T$,
$h_t^\T$ and ${\bm1}^\T$ are the transposes of $h_t$ and ${\bm1}$.
An alternative form of \eqref{e_Da} is $D_t=(U_t+U_t^{\xy\T})/2$, where $U_t=h_t{\bm1}^\T-H_t$, and its elementwise form is
\eq{
\label{e_da}%+
d_{ij}(t)=\frac{1}{2}\xy(h_{ii}(t)+h_{jj}(t)-h_{ij}(t)-h_{ji}(t)),\quad i,j\in V(G),
}
where $H_t=(h_{ij}(t))$ and $D_t=(d_{ij}(t)).$ This is a standard transformation used to obtain a distance from a proximity measure (cf.\ the inverse covariance mapping in~\cite[Section\:5.2]{DezaLaurent97} and the cosine law in~\cite{Critchley88}).

\medskip
In the rest of this section, we present several known facts (lemmas) which will be of use in what follows, one simple example, and two remarks.

%The first lemma follows from Theorem~6 in~\cite{Che11AAM}.

\begin{lemma}[\!\!\cite{Che11AAM}]
\label{l_di}
For any connected $G,$ if ${R_t=(r_{ij}(t))}$ exists$,$ then the matrix $D_t=(d_{ij}(t))$ defined by~\eqref{e_Rt}--\eqref{e_Da} determines a graph-geodetic distance $d_t(i,j)=d_{ij}(t)$ on~$V(G).$
\end{lemma}

This enables one to give the following definition.

\begin{defin}
\label{d_walkD}
{\rm
For a connected graph $G,$ the \emph{walk distances\/} on~$V(G)$ are the functions $d_t(i,j)\!:V(G)\!\times\! V(G)\to\R$ and the functions, $d^\W_t(i,j),$ positively proportional to them$,$ where $d_t(i,j)=d_{ij}(t)$ %$i,j\in V(G)$
and $D_t=(d_{ij}(t))$ is defined by~\eqref{e_Rt}--\eqref{e_Da}.
}
\end{defin}
%We call the distances described in Lemma~\ref{l_di} and their positive multipliers the \emph{walk distances\/} on~$V(G).$

\begin{example}
\label{ex_1}
{\em
For the multigraph $G$\/ shown in Fig.\:\ref{p_3},
\begin{figure}[th] %[p]
% \hspace{-1.5em}\input GRAPH.LP
\begin{center}
\bigskip
\includegraphics[width=2in]{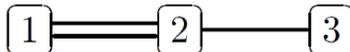}
\end{center}
\Up{.9}
\caption{A multigraph $G$ on 3 vertices.\label{p_3}}
\end{figure}

\noindent the weighted adjacency matrix is
$$
A=\left[\begin{array}{rrr}
0&2&0\\
2&0&1\\
0&1&0\\
\end{array}\right],
$$
the matrix $R_\frac13$ of the walk weights of $\frac13G$ exists and has the form
$$
R_\frac13
=R\bigl(\tfrac13G\bigr)
=\bigl(r_{ij}\bigl(\tfrac13\bigr)\bigr)
=\frac14\left[\begin{array}{rrr}
8&6&2\\
6&9&3\\
2&3&5\\
\end{array}\right],
$$
and the computation \eqref{e_Da} of the walk distances $d_t(i,j)$ with parameter $t=\frac13$ yields %matrix $D_\frac13$ has the form %of walk distances in $G$ with parameter $\,t=1/3\,$ is %\eqref{e_da} yields
$$
D_\frac13
=\bigl(d_{ij}\bigl(\tfrac13\bigr)\bigr)
=\frac12\left[\begin{array}{rrr}
    0&\ln2&\ln10\\
\ln 2&   0&\ln 5\\
\ln10&\ln5&    0\\
\end{array}\right]
\approx\left[\begin{array}{rrr}
   0&0.35&1.15\\
0.35&   0&0.80\\
1.15&0.80&   0\\
\end{array}\right].
$$
Since the walk distances are graph-geodetic (Lemma\:\ref{l_di}) and all paths from $1$ to $3$ visit $2$, $d_\frac13(1,2)+d_\frac13(2,3)=d_\frac13(1,3)$ holds.
}
\end{example}

\medskip
Regarding the existence of $R_t,$ since for a connected graph $A$ is irreducible, the Perron-Frobenius theory of nonnegative matrices provides the following result (cf.~\cite[Theorem~4]{Thompson58}).
\begin{lemma}
\label{l_finite}
For any weighted adjacency matrix $A$ of a connected graph $G,$ the series
$R_t=\sum_{k=0}^\infty(tA)^k$ with $t>0$ converges to $(I-tA)^{-1}$ if and only if\/ $t<\rho^{-1},$ where $\rho=\rho(A)$ is the spectral radius of\/~$A.$ Moreover$,$ $\rho$ is an eigenvalue of~$A;$ as such $\rho$ has multiplicity~$1$ and a positive eigenvector.
\end{lemma}

Observe that for the graph $G$ of Example~\ref{ex_1}, $\rho=\sqrt{5},$ so $\frac13=t<\rho^{-1}$ is satisfied.
%Eigenvalue $\rho=\rho(A)\xy$ is called the \emph{Perron root\/} of~$A.$
%If $x$ is an eigenvector of $A$ associated with~$\rho,$ then the probability vector $p=x/\|x\|_1$ is called the \emph{Perron vector of}~$A.$
%
\begin{lemma}
\label{l_expD}
For any vertices $\,i,j\in V(G)$ and $\,0<t<\rho^{-1},$
\eqss{
\label{e_d1}%+
d_t(i,j)
=-\ln\biggl(\xz\frac{r_{ij}(t)}{\sqrt{r_{ii}(t)\,r_{\!jj}(t)}}\xz\biggr).
}
\end{lemma}

Lemma~\ref{l_expD} is a corollary of \eqref{e_Ha}, \eqref{e_Da}, and Lemma\:\ref{l_finite}. %(cf.\ Eq.~(11) in~\cite{Che11AAM}).

\medskip
On the basis of Lemma\:\ref{l_expD}, the walk distances can be given the following short definition: $d_t(i,j)=-\ln r'_{ij}(t),$ where $r'_{ij}(t)=\frac{r_{ij}(t)}{\sqrt{r_{ii}(t)\,r_{\!jj}(t)}}$ and $R_t=(r_{ij}(t))_{n\times n}$ is defined by~\eqref{e_Rt}.

\begin{remark}
{\rm
Consider another transformation of the correlation-like index $r'_{ij}(t)=\frac{r_{ij}(t)}{\sqrt{r_{ii}(t)\,r_{\!jj}(t)}}$:
\eq{
\label{e_d2}
d'_t(i,j)
=1-\frac{r_{ij}(t)}{\sqrt{r_{ii}(t)\,r_{\!jj}(t)}}. %,\quad i,j\in V(G).
}
%Does $d'_t(i,j)$ have any interesting properties?
Is $d'_t(i,j)$ a metric? It follows from Definition\:\ref{d_walkD}, \eqref{e_d1}, and \eqref{e_d2} that for any walk distance $d^\W_t(i,j),$ there exists $\lambda>0$ such that
\eq{
\label{e_Sho}
d'_t(i,j)
=1-e^{-\lambda d^{\mathrm{w}}_t\xz(i,\xy j)}.
}

Eq.\:\eqref{e_Sho} is the Schoenberg transform~\cite{Schoenberg35,Schoenberg38} (see also \cite[Section\:9.1]{DezaLaurent97} and \cite{Laurent98,Bavaud11}). As mentioned in \cite{DezaDeza09}, an arbitrary function $\td(i,j)$ is the result of the Schoenberg transform of some metric if and only if $\td(i,j)$ is a \emph{P-metric\/}, i.e., a metric with values in $[0,\,1]$ that satisfies the \emph{correlation triangle inequality\/}
\eq{
1-\td(i,k)\ge(1-\td(i,j))(1-\td(j,k)),
}
which can be rewritten as $\td(i,k)\le \td(i,j)+\td(j,k)-\td(i,j)\xy\xy\td(j,k).$

This fact implies that \eqref{e_d2} defines a P-metric. It is easily seen that the correlation triangle inequality for $d'_t(i,j)$ reduces to the transition inequality~\eqref{e_te}; obviously, it can be given a probabilistic interpretation.

For the graph $G$ of Example\:\ref{ex_1}, the P-metric $d'_t(i,j)$ with $t=\frac13$ is given by the matrix
$$
D'_\frac13
=\bigl(d'_{ij}\bigl(\tfrac13\bigr)\bigr)
=\left[\begin{array}{rrr}
           0&1-\sqrt{0.5}&1-\sqrt{0.1}\\
1-\sqrt{0.5}&           0&1-\sqrt{0.2}\\
1-\sqrt{0.1}&1-\sqrt{0.2}&           0\\
%          0&1-2^{-0.5}&1-10^{-0.5}\\
%1- 2^{-0.5}&         0&1- 5^{-0.5}\\
%1-10^{-0.5}&1-5^{-0.5}&          0\\
\end{array}\right]
\approx\left[\begin{array}{rrr}
   0&0.29&0.68\\
0.29&   0&0.55\\
0.68&0.55&   0\\
\end{array}\right].
$$
}
\end{remark}

\begin{remark}
{\rm
It can be noted that the Nei standard genetic distance \cite{Nei72} and the Jiang-Conrath semantic distance \cite{JiangConrath97} have a form similar to~\eqref{e_d1}.
Moreover, the transformation $\,-\ln(r(i,j))$ where $r(i,j)$ is a similarity measure between objects $i$ and $j$ was used in the construction of the Bhattacharyya distance between probability distributions~\cite{Bhattacharyya43} and the Tomiuk-Loeschcke genetic distance \cite{TomiukLoeschcke91} (see also the Leacock-Chodorow similarity \cite{LeacockChodorow98} and the Resnik similarity~\cite{Resnik95}). These and other distances and similarities are surveyed in~\cite{DezaDeza09}.
}
\end{remark}

%\Up{.3}
\section{An interpretation of the walk distances} %in graphs}
\label{s_interp}

%Let $n>2.$
For a fixed $t\!:0<t<\rho^{-1},$ where $\rho=\rho(A)$ let us use the notation
\eq{
\label{e_B}
B=I-tA.
}

Assume that $i$ and $j\ne i$ are also fixed and that $i+j$ is even; otherwise this can be achieved by renumbering the vertices. Hence, using \eqref{e_Rt}--\eqref{e_da}, 
the positivity of $R_t=(I-tA)^{-1},$ 
and the determinant representation of the inverse matrix
we obtain
\eq{
\label{e_1}
d_t(i,j)
=0.5(\ln\det B_\ii+\ln\det B_\jj-\ln\det B_\ij-\ln\det B_\ji),
}
where $B_\ij$ is $B$ with row $i$ and column $j$ removed. %and $\th$ is a scaling factor.

\subsection{Logarithms of the cofactors: expressions in terms of circuits} % of $I-tA$
\label{ss_interp1}

To obtain an interpretation of the right-hand side of \eqref{e_1}, we need the following remarkable result due to Kasteleyn.

\begin{lemma}[Kasteleyn \cite{Kasteleyn67}] %[Section\:III.D]
\label{l_K}
For a digraph\/ $\G$ with a weighted adjacency matrix $\tilde A,$
\eqss{
\label{e_Bi1}
\det(I-\tilde A)
&=&\exp\left(-\sum_{c\xy\in\xy\CC}\frac{w(c)}{\mu(c)}\right)\\
\label{e_Bi2}
&=&\prod_{c\xy\in\xy\CC_1}(1-w(c)),
}
where $\CC$ and $\CC_1$ are the sets of all circuits and of all non-periodic circuits in $\G,$ $w(c)$ and $\mu(c)$ being the weight and the multiplicity of the circuit~$c.$
\end{lemma}

The representation \eqref{e_Bi1} was obtained by considering the generating function of walks in~$\G.$
%The sum $\sum_{c\xy\in\xy\CC}\frac{w(c)}{\mu(c)}$ can be treated as a formal counting series in the %counting variables corresponding to abstract arc ``weights''
Basically, the sum $\sum_{c\xy\in\xy\CC}\frac{w(c)}{\mu(c)}$ is a formal counting series in abstract \emph{weight variables\/} (cf.\:\cite[p.\:19]{Riordan58}). However, as soon as the weights are real and thus the generating function is a function in real counting variables, %rather than a formal counting series,
the issue of convergence arises. Since \eqref{e_Bi1} is based on the power expansion $\xy-\ln(I-\tilde A)=\sum_{k_{}=1}^\infty k^{-1}\tilde A^k,$ a necessary condition of its validity in the real-valued setting is $\rho(\tilde A)<1.$

When the arc weights are nonnegative, the same condition is sufficient. However, if %some arc weights are negative and
some vertices $i$ and $j$ are connected by parallel $i\too j$ arcs carrying weights of different signs, then the problem of conditional convergence arises. Namely, if the absolute values of such weights are large enough, then, even though $\rho(\tilde A)<1,$ by choosing the order of summands in the right-hand side of \eqref{e_Bi1}, the sum can be made divergent or equal to any given number.

To preserve \eqref{e_Bi1} in the latter case, the order of summands must be adjusted with an arbitrary order of items in $\sum_{k_{}=1}^\infty k^{-1}\tilde A^k.$ Hence it suffices to rewrite \eqref{e_Bi1} in the form
\eq{
\label{e_Bi1a}
\det(I-\tilde A)
=\exp\left(-\sum_{k=1}^\infty\sum_{c\xy\in\xy\CC_k}\frac{w(c)}{\mu(c)}\right),
}
where $\CC_k$ is the set of all circuits that involve $k$ arcs in $\G.$

\medskip
Lemma~\ref{l_K} is also applicable to undirected graph. To verify this, it is sufficient to replace an arbitrary undirected graph $G$ with its \emph{directed version}, i.e., the digraph obtained from $G$ by replacing every edge by two opposite arcs carrying the weight of that edge.

\medskip
Since by \eqref{e_B}, $B_\ii=I-(tA)_\ii,$ Lemma\:\ref{l_K} can be used to evaluate $\,\ln\det B_\ii.$ Let $G_\i$ ($G_\ij$) be $G$ with vertex $i$ (vertices $i$ and $j$) and all edges incident to $i$ ($i$ and $j$) removed.

\begin{corol}
\label{c_Bii}$\,$

\Up{2.2}
\eqs*{
-\ln\det B_\ii
=\sum_{c\xy\in\xy\CC^{\i}}              \frac{w(c)}{\mu(c)}
=\sum_{c\xy\in\xy\CC^{\ij}\cup\CC^{j\i}}\frac{w(c)}{\mu(c)},
}

\Up{1}\noindent
where
\begin{itemize}
\Up{.50}
\item $\CC^{\i} $ is the set of circuits                               in\/ $tG_\i,$
\Up{.75}
\item $\CC^{\ij}$ is the set of circuits                               in\/ $tG_\ij,$
\Up{.75}
\item $\CC^{j\i}$ is the set of circuits visiting $j,$ but not $i$ in\/ $tG,$
\Up{.50}
\end{itemize}
$w(c)$ and $\mu(c)$ being the weight and the multiplicity of\/~$c.$
\end{corol}

\proof
By assumption, $0<t<\rho^{-1}(A);$ $\,B_\ii=I-tA_\ii.$ %therefore, $\rho(tA)<1.$
Since $A$ is irreducible, $\rho(tA_\ii)<\rho(tA)<1$ \cite[Ch.\:III, \S\:3.4]{Gantmacher59-2}. %\cite[item\:5.5.6 of Part\:II]{MarcusMinc64}.
Moreover, the edge weights in $G$ are positive by assumption.  Therefore, the expansion \eqref{e_Bi1} holds for $B_\ii,$ which yields the desired statement.
\qed\medskip

%To ensure convergence in the expression Lemma~\ref{l_K} provides for $\,\ln\det B_\ij,$ an appropriate transformation of $B_\ij$ is necessary.
To interpret \eqref{e_1}, we also need an expansion of $\,\ln\det B_\ij$ ($j\ne i$). Convergence in such an expansion provided by Lemma~\ref{l_K} can be achieved by applying a suitable linear transformation of~$B_\ij.$
%Convergence in the expression provided by Lemma~\ref{l_K} for $\,\ln\det B_\ij$ ($j\ne i$) can be achieved by applying a suitable linear transformation of~$B_\ij.$

For the fixed $i$ and $j\ne i,$ consider the matrix
\eq{
\label{e_Tij}%
T_{ij}=I(j,i)\_\ij,
}
where $I(j,i)$ differs from $I_{n\times n}$ by the $ji$-entry: $I(j,i)_{ji}=-1.$ 

%\smallskip
The reader can easily construct examples of $T_{ij}$ and verify the following properties.

%$T_{ij}$~has the following properties.
\begin{lemma}
\label{l_Tij}$\,$

$1.$ The columns of  $\,T_{ij}$ form an orthonormal set$,$ i.e.$,$ $T_{ij}$ is orthogonal\/$:$ $T_{ij}^{\xy\T}\xy T_{ij}=I$.

$2.$ If $\,i+j$ is even\/ $(\xz$as assumed$\xy),$ then $\,\det T_{ij}=1.$

$3.$  $\,T_{ij}^{\xy\T}=T_{ji}.$

$4.$ For any $M_{n\times n},\,$ $M_\ij\xy T_{ij}^{-1}$ is obtained from $M$ by\/$:$
{\rm  (i)}\:deleting row $i,$\/
{\rm (ii)}\:multiplying column $i$ by~$-1,$ and\/
{\rm(iii)}\:moving it into the position of column~$j.$
\end{lemma}

The proof of Lemma\:\ref{l_Tij} is straightforward.

\begin{corol}
\label{c_Tij}
$1.$ $I_\ij\xy T_{ij}^{-1}$ is obtained from $I_{(n-1)\times (n-1)}$ by replacing the $kk$-entry with\/ $0,$ where
\eq{
\label{e_kk}
\,k=\begin{cases}j,&j<i,\\j-1,&j>i.\end{cases}
}
$2.$ $I_\ij\xy T_{ij}^{-1}I_\ij=I_\ij,$ i.e.$,$ $T_{ij}^{-1}$ is a g-inverse {\rm\cite{RaoMitra71}} of\/~$I_\ij.$
\end{corol}

Since $\det T_{ji}=1$ (Lemma~\ref{l_Tij}), we have
\eq{
\label{e_detB}
\det B_\ij
=\det (B_\ij\xy T_{ji}).
}

Now we apply Kasteleyn's Lemma\:\ref{l_K} to $B_\ij\xy T_{ji}$ by considering a (multi)digraph $\G$ whose weighted adjacency matrix is
\eq{
\label{e_AA}
\AA=I-B_\ij T_{ji},
}
where $B$ is defined by \eqref{e_B}. Namely, Lemma\:\ref{l_K} in the form \eqref{e_Bi1a} along with \eqref{e_detB} yield

\begin{lemma}
\label{l_1}
\eq{
\label{e_detBij}
-\ln\det B_\ij
=\sum_{k=1}^\infty\sum_{c\xy\in\xy\CC'_k}\frac{w(c)}{\mu(c)},
}
where $\CC'_k$ is the set of all circuits that involve $k$ arcs in a digraph\/ $\G$ whose weighted adjacency matrix is~$\AA,$
while $w(c)$ and $\mu(c)$ are the weight and the multiplicity of the circuit~$c.$
\end{lemma}

As well as \eqref{e_Bi1}, \eqref{e_detBij} is applicable to the case of formal counting series. However, in \eqref{e_B}, $t$ is a real weight variable. In this case, a necessary and sufficient condition of the convergence in \eqref{e_detBij} is $\rho(\AA)<1.$
%It should be noted that $t<\rho^{-1}$ does not generally imply this.
%Recall that $t$ in \eqref{e_B} is a real variable. In this case, a necessary and sufficient condition of the convergence

%It remains to
Let us clarify the relation of $\G$ and its circuits with $G$ and its topology.
This is done in the following section. %by the following statements.

\subsection{The walk distances: An expression in terms of walks}
\label{ss_interp2}

To elucidate the structure of the digraph $\G$ introduced in Lemma\:\ref{l_1}, an algorithmic description of the matrix $\AA$ is useful.

\begin{lemma}
\label{l_Adj1}
$\AA$ can be obtained from $tA$ by\/$:$
%{\rm  (i)}\:
replacing\/ $ta_{ji}$ with\/ $ta_{ji}-1,$
deleting row\/ $i,$\/
%{\rm (ii)}\:
multiplying column\/ $i$\/ by\/~$-1,$ and\/
%{\rm(iii)}\:
moving it into the position of column~$j.$
\end{lemma}

\proof
By \eqref{e_AA}, items\:1 and\:3 of Lemma\:\ref{l_Tij}, \eqref{e_Tij}, and \eqref{e_B} we have
$$
\AA
=(T_{ij}  -   B_\ij)\xy T_{ij}^{-1}
=(I(j,i)-I+tA)\_\ij \xy T_{ij}^{-1}.
$$
Now the result follows from item\:4 of Lemma\:\ref{l_Tij}.
\qed\bigskip

Let us reformulate Lemma\:\ref{l_Adj1} in terms of $G$ and~$\G.$
Recall that a digraph is the \emph{directed version\/} of a graph if it is obtained by replacing every edge in the graph by two opposite arcs carrying the weight of that edge. %they share the same weighted adjacency matrix.

\begin{corol}
\label{c_Adj1}
A digraph\/ $\G$ with weighted adjacency matrix $\AA$ can be obtained from $tG$ by\/$:$
\begin{itemize}
\Up{.75}
\item taking the directed version of the restriction of\/ $tG$ to $V(G)\sm\{i,j\}\xy$ and
\Up{.75}
\item adding a vertex $ij$ with\/$:$ two loops of weights\/ $1$ and\/ $-ta_{ji}$ $(\xz$negative\/\footnote{If $a_{ji}=0,$ then this loop is omitted.}$),$
                                weights\/ $ta_{jm}$ of outgoing arcs\/$,$ and
                                weights\/ $-ta_{mi}$ of incoming arcs$,$ where\/ $m\in V(G)\sm\{i,j\}.$
\end{itemize}
\Up{.75}
Vertex\/ $ij$ is represented in $\AA$ by row and column $k,$ where $k$ is given by~\eqref{e_kk}.
\end{corol}

In what follows, $\G$ denotes the digraph defined in Corollary\:\ref{c_Adj1}.
The \emph{jump\/} in $\G$ is the loop of weight $1$ at~$ij.$
The walk in $\G$ that consists of one jump is called the \emph{jump walk\/} (at~$ij$).

To interpret $\ln\det B_\ij$ in terms of $G,$ we need the following notation.

\begin{defin}
\label{d_wwj}
{\em
A \emph{walk with $i,j$ jumps in $G$\/} is any walk in the graph $G'$ obtained from $G$ by attaching two additional loops of weight $1${\rm:} one adjacent to vertex $i$ and one adjacent to~$j.$ These loops are called \emph{jumps.}
A~walk with $i,j$ jumps (in $G$) only consisting of one jump is called a \emph{jump walk\/} (at $i$ or~$j$).
}
\end{defin}

\begin{defin}
\label{d_owwj}
{\em
A \emph{$j\too i$ alternating walk with jumps\/} is any $j\too i$ walk $\w$ with $j,i$ jumps such that
%a $j\ldots j$ subwalk of\/ $\w$ may contain no occurrence of\/ $i$ only if all edges it contains are jumps
$($a$)$~any $j\ldots j$ subwalk of\/ $\w$ either visits $i$ or contains no edges except for jumps and
$($b$)$~any $i\ldots i$ subwalk of\/ $\w$ either visits $j$ or contains no edges except for jumps.
%similarly for the appearance of $j$ on the $i\ldots i$ subwalks.
%a $i\ldots i$ subwalk of $\w$ may contain no occurrences of $j$ only if all edges it contains are jumps

A \emph{$j\too i\too j$ alternating walk with jumps\/} is defined similarly: the only difference is that the endpoint of such a walk is~$j.$
}
\end{defin}

%Furthermore, it can be observed
To introduce some additional notation, observe that any $j\too i$ alternating walk $\w$ with jumps can be uniquely partitioned into a sequence of subwalks $(\w_1\cdc \w_t)$ such that every two neighboring subwalks share one terminal vertex and each $\w_k$ is a jump walk or is a $j\too i$ or an $i\too j$ hitting walk without jumps. %has no edges except for jumps.
For every $k\in\{1\cdc t\},$ consider the set ${p_k=\{\w_k,\tilde\w_k\},}$ where $\tilde\w_k$ is \underline{either}
$\w_k$ written from end to beginning (reversed\footnote{Cf. ``dihedral equivalence'' in \cite{HararySchwenk79}.}) when $\w_k$ is a hitting walk without jumps, \underline{or}
a jump walk at $i$ ($j$) when $\w_k$ is a jump walk at $j$ (resp.,~$i$).
%records a series of $m$ jumps at $i$ whenever $\w_k$ records a series of $m$ jumps at $j,$ and vice versa.
The sequence $p(\w)=(p_1\cdc p_t)$ will be called the \emph{route partition of\/~$\w.$} We say that two $j\too i$ alternating walks with jumps, $\w$ and $\w',$ are \emph{equipartite\/} if the route partition of $\w'$ can be obtained from that of $\w$ by a cyclic shift. Finally, any equivalence class of equipartite $j\too i$ alternating walks with jumps will be called an \emph{alternating $j\too i$ route with jumps.} If $\rr$ is such a route, then its \emph{length\/} and \emph{weight\/} are defined as the common length and weight of all walks with jumps it includes, respectively. If a route partition $p(\w)=(p_1\cdc p_t)$ has period (the length of the elementary repeating part)~$y,$ then the \emph{multiplicity\/} of the alternating $j\too i$ route with jumps that corresponds to $p(\w)$ is defined to be~$t/y.$

Completely the same construction can be applied to define \emph{alternating $j\too i\too j$ route with jumps} (starting with the above definition of a $j\too i\too j$ alternating walk with jumps). A~notable difference is that there are alternating $j\too i\too j$ routes with jumps that do not visit~$i$: these consist of jumps at~$j.$ %only involve jumps at~$j.$
The weight of such a route with jumps is $1$ and its multiplicity is the number of jumps.

\begin{lemma}
\label{l_G-G}
There is a one-to-one correspondence between the set of circuits in $\G$ that contain vertex $ij$ and have odd $(\xz$even$)$ numbers of negatively weighted arcs and the set of alternating $j\too i$ routes $(\xz$alternating $j\too i\too j$ routes$)$ with jumps in~$G.$ The circuit in $\G$ and route with jumps in $G$ that correspond to each other have the same length\/$,$ weight\/$,$ and multiplicity.
\end{lemma}

\proof
Every circuit containing vertex $ij$ in $\G$ can be uniquely represented by a cyclic sequence\footnote{A \emph{cyclic sequence\/} is a set $X=\{x_1\cdc x_N\}$ with the relation ``\emph{next}'' $\eta=\{(x_2,x_1)\cdc(x_N,x_{N-1}),$ $(x_1,x_N)\}.$} of walks each of which either is an\/ $ij\too ij$ walk including %no intermediate occurrence of\/~$ij,$
exactly one negatively weighted arc, or is the jump walk at~$ij.$ %records a series of jumps at~$ij.$
Such a cyclic sequence uniquely determines an alternating $j\too i$ or $j\too i\too j$ route with jumps in~$G$ (if the number of negatively weighted arcs involved in the circuit is odd or even, respectively). %depending on the parity of

On the other hand, every set $p_k=\{\w_k,\tilde\w_k\}$ involved in %such an alternating route with jumps
                                                                    an alternating $j\too i$ or $j\too i\too j$ route with jumps
in $G$ uniquely determines either an\/ $ij\too ij$ walk containing exactly one negatively weighted arc, %no intermediate occurrence of\/~$ij$
or %a series of jumps at\/~$ij$ in~$\G.$
the jump walk at\/~$ij$ in~$\G.$ Thereby, every alternating route with jumps under consideration uniquely determines a circuit in~$\G.$ Furthermore, the two correspondences described above are inverse to each other. Thus, these reduce to a one-to-one correspondence.

Finally, it is easily seen that the corresponding circuits and alternating routes with jumps share the same length$,$ weight$,$ and multiplicity.
\qed%\medskip

\begin{remark}
{\rm
It can be noted that the multiplicity of an alternating $j\too i$ route with jumps in~$G$ can only be odd.
}
\end{remark}

Both circuits and alternating routes will be called \emph{figures}. Lemmas\;\ref{l_1} and\;\ref{l_G-G} enable one to express $\ln\det B_\ij$ in terms of figures in $tG$ and $tG_\ij.$

\begin{lemma}
\label{l_Bij}
\eqs*{
%\label{e_detBij1}
-\ln\det B_\ij\;
=\;\sum_{k=1}^\infty\;\sum_{c\xy\in\xy(\CC^{\ij}\cup\CC^{j\to i\to j}\cup\CC^{j\to i})\xy\cap\,\CC_k}\!\!\!
     %{c\xy\in\xy\CC^{\ij}\cup\CC^{i\tot j}\cup\CC^{i\mi j}}
 (-1)^{\ze(c)}\xy\frac{w(c)}{\mu(c)},
}
where %$\CC^{[\cdots\hspace*{-.01em}]}_k$ are the subsets comprising figures that involve $k$ arcs of\/$:$
\begin{itemize}
%\Up{.75}
\item $\CC^{\ij}$          is the set of                 circuits            in\/ $tG_\ij,$
\Up{.75}
\item %$\CC^{i\tot j}$ is the set of      $i\tott j$ circuits with jumps in\/ $tG,$
       $\CC^{j\to i\to j}$ is the set of alternating\/ $j\too i\too j$ routes with jumps in\/ $tG,$
\Up{.75}
\item %$\CC^{i\mi  j}$ is the set of open $i\mii  j$ routes   with jumps in\/ $tG,$
       $\CC^{j\to i}$      is the set of alternating\/ $j\too i$ routes with jumps in\/ $tG,$
\Up{.75}
\item  $\CC_k$      is the set of figures $(\xz$in\/ $tG$ or $tG_\ij)$ that involve $k$ arcs$,$
\Up{.6} %{1.25}
\end{itemize} %$\CC^{\ij}_k,$ $\CC^{j\to i\to j}_k,$ and $\CC^{j\to i}_k$
$$
\ze(c)=\begin{cases}0,&c\in\CC^{\ij}\cup\CC^{j\to i\to j},\\ %{i\tot j}
                    1,&c\in\CC^{j\to i}, %{i\mi j}
       \end{cases}
$$
while $w(c)$ and\/ $\mu(c)$ are the weight and the multiplicity of\/~$c.$
%$\,\CC^{\ij}=\cup_{k=1}^\infty\CC^{\ij}_k$ and similarly for the other sets of circuits and alternating routes.
\end{lemma}

Similarly, we can express $\ln\det B_\ji$ in terms of the sets $\CC^{\ij},$ $\CC^{i\to j\to i},$ and $\CC^{i\to j}.$ There exist natural bijections between $\CC^{j\to i\to j}$ and $\CC^{i\to j\to i}$ and between $\CC^{j\to i}$ and $\CC^{i\to j}.$ Namely, to obtain an element of\/ $\CC^{i\to j\to i}$ from $c\in\CC^{j\to i\to j}$ (or an element of\/ $\CC^{i\to j}$ from $c\in\CC^{j\to i}$), it suffices to reverse all $j\too i$ and $i\too j$ hitting walks without jumps in $c$ and to replace every %series of jumps %of the same length,
jump walk at $j$ with the jump walk at $i$ and vice versa.
%A bijection between $\CC^{j\to i}$ and $\CC^{i\to j}$ is established similarly.

On the other hand, the sets $\CC^{i\tot j}\stackrel{\rm def}{=}\CC^{j\to i\to j}\cup\CC^{i\to j\to i}$
                        and $\CC^{i\mi  j}\stackrel{\rm def}{=}\CC^{j\to i}     \cup\CC^{i\to j}$ also make sense.
Specifically, they are useful for expressing $d_t(i,j).$ %, see Theorem\:\ref{t_dij}.
Such an expression is the main result of this paper. It follows by combining \eqref{e_1}, Corollary\:\ref{c_Bii}, and Lemma\:\ref{l_Bij}.
%Combining \eqref{e_1}, Corollary\:\ref{c_Bii}, and Lemma\:\ref{l_Bij} we obtain the main result of this paper.
\begin{thm}
\label{t_dij}
$$d_t(i,j)\;
=\;\frac12\,\sum_{k=1}^\infty\;\sum_{c\xy\in\xy(\CC^{i\j}\cup\CC^{\i j}\cup\CC^{i\tot j}\cup\CC^{i\mi j})\xy\cap\,\CC_k}\!\!\!
            %{c\xy\in\xy\CC^{i\j}\cup\CC^{\i j}\cup\CC^{j\to i\to j}\cup\CC^{j\to i}}
  (-1)^{\ze(c)}\xy\frac{w(c)}{\mu(c)},$$
where the \underline{sets} of figures in $tG$ are denoted by\/$:$
%$\CC^{[\cdots\hspace*{-.01em}]}_k$ are the sets of circuits or alternating routes that involve $k$ arcs\/$:$
\begin{itemize}
\Up{.75}
\item $\CC^{i \j}\!:$     of                 circuits visiting\/ $i,$ but not\/ $j,$
\Up{.75}
\item $\CC^{\i j}\!:$     of                 circuits visiting\/ $j,$ but not\/ $i,$
\Up{.75}
\item  $\CC^{i\tot j}\!:$ of alternating\/ $j\too i\too j$ and\/ $i\too j\too i$ routes with jumps\/$,$
      %$\CC^{j\to i\to j}$ of alternating\/ $j\too i\too j$ routes with jumps in\/ $tG,$
      %  $i\tott j$ circuits with jumps in\/ $tG,$
\Up{.75}
\item $\CC^{i\mi  j}\!:$  of alternating\/ $j\too i$ and\/ $i\too j$ routes with jumps\/$,$
     %$\CC^{j\to i}$ is the set of alternating\/ $j\too i$ routes with jumps in\/ $tG,$
     %open $i\mii  j$ routes   with jumps in\/ $tG,$
\Up{.75}
\item  $\CC_k\!:$         of figures that involve $k$ arcs\,$;$
\Up{0.6} %\Up{1.25}
\end{itemize}
%\Up{0.6} %\Up{1.25}
%\Up{.5}
$$
\ze(c)=\begin{cases}0,&c\in\CC^{i\tot j},\\ %{j\to i\to j}
                    1,&c\in\CC^{i\j}\cup\CC^{\i j}\cup\CC^{i\mi j}, %{j\to i}
       \end{cases}
$$
while $w(c)$ and\/ $\mu(c)$ are the weight and the multiplicity of\/~$c.$ %and $\th$ is a scaling factor.
\end{thm}

In more general terms, Theorem\:\ref{t_dij} can be interpreted as follows. The walk distance between $i$ and $j$ is reduced by $j\too i$ and $i\too j$ walks (see $\CC^{i\mi j}$), connections of $i$ with other vertices avoiding $j$ ($\CC^{i\j}$), and connections of $j$ avoiding $i$ ($\CC^{\i j}$).
The set $\CC^{i\tot j}$ supplies all positive terms in the expansion of $d_t(i,j).$ It comprises constantly jumping walks along with closed walks involving $i$ and $j$ whose positive weights compensate the negative overweight of $j\too i$ and $i\too j$ routes with extra jumps. 

Note that Theorem\:\ref{t_dij} supports the observation in the Introduction that the high centrality of $i$ and $j$ reduces, ceteris paribus, the walk distance between them. Indeed, the elements of $\CC^{i\j}\cup\CC^{\i j}$ which account for the centrality of $i$ and $j$ make a negative contribution to the distance.

The following example may provide some additional insight into Theorem\:\ref{t_dij}.

\begin{example}
\label{ex_2}
{\em
For the graph $G$ of Example\:\ref{ex_1}, let us approximate $d_\frac13(1,3)=\frac12\ln10\approx1.15\xy$ using Theorem\:\ref{t_dij}. Due to~\eqref{e_AA},
$\AA=\dfrac13\left[\begin{array}{rr}0&-2\\1& 3\\\end{array}\right].$
As $\rho(\AA)=2/3<1,$ convergence holds in \eqref{e_detBij} %and~\eqref{e_detBij1}
and thus in Theorem\:\ref{t_dij}. The leading terms of the expansion Theorem\:\ref{t_dij} provides for $d_\frac13(1,3)$ are presented in Table\:\ref{ta_ex}.
%For every set $\CC^{[\cdots\hspace*{-.01em}]}$ in Theorem\:\ref{t_dij} by $\CC^{[\cdots\hspace*{-.01em}]}_k$ we denote $\CC^{[\cdots\hspace*{-.01em}]}\cup\CC_k.$
In this table, $\dfrac k\mu(v_0\cdots v_m)$ is the denotation of a collection of figures where each figure has multiplicity $\mu$ and contains some walk (or walk with jumps) whose sequence of vertices is $v_0\cdc v_m$; $k$ is the cardinality of the collection.
%$v_0\cdc v_m$ is the sequence of vertices of some walk (or walk with jumps), $k$ is the number of figures that contain a walk (with jumps) having this sequence of vertices, and $\mu$ is the multiplicity of each such a figure.
If $\mu=1,$ then $\mu$ is omitted; if $\mu=k=1,$ then $\mu$ and $k$ are omitted.

\begin{table}[ht]
\Up{.7}
$$
\begin{array}{c||c|c|c|}                                                                                       %\hline
\cap &\CC^{1 \bar3}\cup\CC^{\bar1 3}&\CC\ak^{1\tot 3\aj}                                  &\CC^{1\mi  3}      \\\hline\hline
\CC_1&\varnothing                   &       (11)\ai,           (33)\ai                    &\varnothing        \\\hline
\CC_2&4(121),(323)                  &\frac12\ai(111),   \frac12(333)                      &2 (123),2 (321)    \\\hline
\CC_3&\varnothing                   &\frac13\ai(1111),  \frac13(3333)                     &2(1123),2(3321)    \\\hline
\CC_4&\!\!\!\begin{array}{c}
      \frac42\aj(12121),6(12121),\\
      \frac12\ak(32323)
      \end{array}\!\!\!            &\begin{array}{c}
                                    \frac14\ai(11111), \frac14(33333),\\
                                    \frac22\ak(12321),(12321),\frac22(32123),(32123)
                                    \end{array}                                           & 2(11123), 2(33321)\\\hline
\CC_5&\varnothing                  &\frac15\ai(111111),\frac15(333333),4(112321),4(332123)&2(111123),2(333321)\\\hline
\end{array}
$$
\caption{The figures forming the leading terms in the expansion of $d_\frac13(1,3)$ in Example\:\ref{ex_2}.\label{ta_ex}}
\end{table}

The first terms of the series Theorem\:\ref{t_dij} provides are:
\eqs*{
d_\frac13(1,3)
&\!=\!&\frac12
 \Biggl[(2\xz\cdot\xz1)
+\left(-\frac49-\frac19+2\xz\cdot\xz\frac12-2\xz\cdot\xz\frac29\right)
+\left(2\xz\cdot\xz\frac13-2\xz\cdot\xz\frac29\right)\\
&\!+\!&\left(-\frac{2+6}{81}-\frac12\xz\cdot\xz\frac1{81}+2\biggl(\frac14+\frac{1+1}{81}\biggr)-2\xz\cdot\xz\frac29\right)
+\left(2\biggl(\frac15+\frac4{81}\biggr)-2\xz\cdot\xz\frac29\right)
+\ldots\Biggr]\\
&\!=\!&\frac{461}{405}+\ldots,
}
where $\frac{461}{405}\approx1.1383.$ 

In the above expression, the sum (with signs) of the weights of figures that involve $k$ edges is $0$ whenever $k$ is even. Thus, the above expansion reduces to
\eqs*{
d_\frac13(1,3)
=\frac12
 \Biggl[(2\xz\cdot\xz1)
+\left(2\xz\cdot\xz\frac13-2\xz\cdot\xz\frac29\right)
+\left(2\biggl(\frac15+\frac4{81}\biggr)-2\xz\cdot\xz\frac29\right)
+\ldots\Biggr].
}

The relative error of this approximation is $1.1\%.$ 
}
\end{example}

\medskip
In some cases, the convergence of such expansions is extremely slow. On the other hand, the meaning of Theorem\:\ref{t_dij} is to clarify the concept of walk distance by representing it as the sum of route/circuit weights rather than to provide an effective algorithm for computing it.
%\bigskip

%\section*{Acknowledgements}
%
%This work was partially supported by the RFBR Grant 09-07-00371 and the RAS Presidium Program ``Development of Network and Logical Control''.
%%The authors are grateful to the referees for their comments.

%\vspace{-.1em}
%\begin{thebibliography}{99}

%\newpage
\sloppy
\emergencystretch=5pt    %
\righthyphenmin=3        %
\hfuzz=1.9pt            %
\tolerance=400           %
\hyphenpenalty50         %
\tolerance200 \hbadness200 %

\bibliographystyle{elsart-num-sort} %{siam}%{gost71sP}
\hfuzz=1.9pt
\bibliography{all2}     %%

\begin{thebibliography}{10}
\expandafter\ifx\csname url\endcsname\relax
  \def\url#1{\texttt{#1}}\fi
\expandafter\ifx\csname urlprefix\endcsname\relax\def\urlprefix{URL }\fi

\bibitem{Bavaud11}
F.~Bavaud, On the {Schoenberg} transformations in data analysis: Theory and
  illustrations, Journal of Classification 28~(3) (2011) 297--314.

\bibitem{Bhattacharyya43}
A.~Bhattacharyya, On a measure of divergence between two statistical
  populations defined by their probability distributions, Bulletin of the
  Calcutta Mathematical Society 35 (1943) 99--109.

\bibitem{Che11DAM}
P.~Chebotarev, A class of graph-geodetic distances generalizing the
  shortest-path and the resistance distances, Discrete Applied Mathematics
  159~(5) (2011) 295--302.

\bibitem{Che11AAM}
P.~Chebotarev, The graph bottleneck identity, Advances in Applied Mathematics
  47~(3) (2011) 403--413.

\bibitem{Che12DAM}
P.~Chebotarev, The walk distances in graphs, Discrete Applied Mathematics, In
  press.
\newline\urlprefix\url{http://dx.doi.org/10.1016/j.dam.2012.02.015}

\bibitem{DezaDeza09}
M.~M. Deza, E.~Deza, Encyclopedia of Distances, Springer, Berlin--Heidelberg,
  2009.

\bibitem{DezaLaurent97}
M.~M. Deza, M.~Laurent, {Geometry of Cuts and Metrics, volume 15 of Algorithms
  and Combinatorics}, Springer, Berlin, 1997.

\bibitem{Critchley88}
F.\;Critchley, On certain linear mappings between inner-product and
  squared-distance matrices, Linear Algebra and its Applications 105 (1988)
  91--107.

\bibitem{Gantmacher59-2}
F.~R. Gantmacher, Applications of the Theory of Matrices, Interscience, New
  York, 1959.

\bibitem{Harary69}
F.~Harary, Graph Theory, Addison-Wesley, Reading, MA, 1969.

\bibitem{HararySchwenk79}
F.~Harary, A.~Schwenk, The spectral approach to determining the number of walks
  in a graph, Pacific Journal of Mathematics 80~(2) (1979) 443--449.

\bibitem{JiangConrath97}
J.~J. Jiang, D.~W. Conrath, Semantic similarity based on corpus statistics and
  lexical taxonomy, in: Proceedings of International Conference on Research in
  Computational Linguistics (ROCLING X), Taiwan, 1997, 15\:pp.

\bibitem{Kasteleyn67}
P.~W. Kasteleyn, Graph theory and crystal physics, in: F.~Harary (ed.), Graph
  Theory and Theoretical Physics, Academic Press, London, 1967, pp. 43--110.

\bibitem{Katz53}
L.~Katz, A new status index derived from sociometric analysis, Psychometrika
  18~(1) (1953) 39--43.

\bibitem{Laurent98}
M.~Laurent, {A connection between positive semidefinite and Euclidean distance
  matrix completion problems}, Linear Algebra and its Applications 273~(1-3)
  (1998) 9--22.

\bibitem{LeacockChodorow98}
C.~Leacock, M.~Chodorow, {Combining local context and WordNet similarity for
  word sense identification}, in: C.~Fellbaum (ed.), WordNet. An electronic
  lexical database, chap.~11, MIT Press, Cambridge, MA, 1998, pp. 265--283.

\bibitem{Nei72}
M.~Nei, Genetic distance between populations, The American Naturalist 106~(949)
  (1972) 283--292.

\bibitem{RaoMitra71}
C.~R. Rao, S.~K. Mitra, Generalized Inverse of Matrices and its Applications,
  Wiley, New York, 1971.

\bibitem{Resnik95}
P.~Resnik, Using information content to evaluate semantic similarity, in:
  Proceedings of the 14th International Joint Conference on Artificial
  Intelligence (IJCAI'95), vol.~1, Morgan Kaufmann Publishers, San Francisco,
  CA, 1995.

\bibitem{Riordan58}
J.~Riordan, An Introduction to Combinatorial Analysis, Wiley, New York, 1958.

\bibitem{Schoenberg35}
I.~J. Schoenberg, Remarks to {M.\:Fr\'echet's} article ``{Sur la d\'efinition
  axiomatique d'une classe d'espaces vectoriels distanci\'es applicables
  vectoriellement sur l'espace de Hilbert}'', Annals of Mathematics 36 (1935)
  724--732.

\bibitem{Schoenberg38}
I.~J. Schoenberg, Metric spaces and positive definite functions, Transactions
  of the American Mathematical Society 44 (1938) 522--536.

\bibitem{Thompson58}
G.~L. Thompson, Lectures on Game Theory, Markov Chains and Related Topics,
  {Monograph SCR}--11, Sandia Corporation, Albuquerque, NM, 1958.

\bibitem{TomiukLoeschcke91}
J.~Tomiuk, V.~Loeschcke, A new measure of genetic identity between populations
  of sexual and asexual species, Evolution 45 (1991) 1685--1694.

\end{thebibliography}
\end{document}